\documentclass[12]{article}
\usepackage{amsfonts}
\usepackage{amsmath}
\usepackage[unicode,colorlinks,plainpages=false]{hyperref}
\usepackage{theorem}

\newtheorem{theorem}{Theorem}
\newtheorem{corollary}{Corollary}
\newtheorem{lemma}{Lemma}

\newcommand{\openbox}{$\begin{array}{c}
\hspace*{-0.55em}\sqcap \hspace*{-0.60em}\\[-0.4em] \hline
\multicolumn{1}{c}{\hspace*{-0.60em}}\\[-0.8em]
\end{array}
$}

\usepackage{enumerate}

\begin{document}

\centerline{{\bf On Special Semigroups Derived From an Arbitrary Semigroup} }

\bigskip
\centerline{Attila Nagy}
\medskip
\centerline{Department of Algebra}
\centerline{Budapest University of Technology and Economics}
\centerline{P.O. Box 91}
\centerline{1521 Budapest}
\centerline{Hungary}
\centerline{e-mail: nagyat@math.bme.hu}

\bigskip

\begin{abstract}
Let $S$ be a semigroup, $\Lambda$ a non-empty set and $P$ a mapping of $\Lambda$ into $S$. The set $S\times \Lambda$ together with the operation
$\circ _P$ defined by $(s, \lambda)\circ _P(t, \mu )=(sP(\lambda)t, \mu )$ form a semigroup which is denoted by $(S, \Lambda , \circ _P)$.
Using this construction, we prove a common connection between the semigroups $S$, $S/\theta$ and $S/\theta ^*=(S/\theta)/(\theta ^*/\theta)$, where $\theta$ and $\theta ^*/\theta$ are the kernels of the right regular representations of $S$ and $S/\theta$, respectively. We also prove an embedding theorem for the semigroup $(S, S/\theta , \circ _p)$, where $S$ is a left equalizer simple semigroup without idempotents, and $P$ maps every $\theta$-class of $S$ into itself.
\end{abstract}

 \section{Introduction}

 \bigskip

Let $S$ be an arbitrary semigroup. It is known that the relation $\theta$ on $S$ defined by $(a, b)\in \theta$ if and only if $xa=xb$ for all $x\in S$ is a congruence on $S$. This congruence is the kernel of the right regular representation $\varphi : a \mapsto \varrho _a$ ($a\in S$) of $S$;
$\varrho _a: s\mapsto sa$ ($s\in S$) is the inner right translations of $S$ defined by $a$. For convenience, the semigroup $\varphi (S)=S/\theta$ is also called the right regular representation of $S$. The $\theta$-class of $S$ containing an element $s\in S$ will be denoted by $[s]_{\theta}$.

Let $\theta ^*$ denote the congruence on the semigroup $S$ for which $\theta\subseteq \theta ^*$ and $\theta ^*/\theta$ is the kernel of the right regular representation on $S/\theta$, where $\theta ^*/\theta$ is defined by $([s]_{\theta}], [t]_{\theta})\in \theta ^*/\theta$ if and only if
$(s, t)\in \theta ^*$ (see Theorem 5.6 of \cite{Howie:sg-11}). It is easy to see that $(a,b)\in \theta ^*$ if and only if $(xa, xb)\in \theta$ for all $x\in S$, that is, $sa=sb$ for all $s\in S^2$ (see also \cite{Nagyleft:sg-6} and \cite{NagyKolibiar:sg-19}). The $\theta ^*$-class of $S$ containing an element $s\in S$ will be denoted by $[s]_{\theta ^*}$.

The right regular representation of semigroups plays an important role in the examination of semigroups. Here we cite some results of papers \cite{Cr1:sg-3}, \cite{Cr2:sg-4} and \cite{Nagy:sg-8}, in which special types of semigroup are characterized by the help of their right regular representation.

A semigroup satisfying the identity $ab=a$ (resp. $ab=b$) is called a left zero (resp. right zero) semigroup. A semigroup is called a left group (resp. right group) if it is a direct product of a group and a left zero (resp. right zero) semigroup.

\begin{lemma}\label{lm1}(\cite{Cr1:sg-3}) A semigroup $S$ is a left group if and only if the right regular representation $S/\theta$ of $S$ is a group.
\end{lemma}

\medskip

A semigroup $S$ is called an $M$-inversive semigroup (\cite{Y:sg-10}) if, for each $a\in S$, there are elements $x, y\in S$ such that $ax$ and $ya$ are middle units of $S$, that is, $caxd=cd$ and $cyad=cd$ is satisfied for all $c, d\in S$.

\begin{lemma}\label{lm2}(\cite{Cr2:sg-4}) A semigroup $S$ is $M$-inversive if and only if the right regular representation $S/\theta$ of $S$ is a right group.
\end{lemma}

\medskip

In \cite{Nagy:sg-8}, a semigroup $S$ is called a left equalizer simple semigroup if, for arbitrary elements $a, b\in S$, the assumption $x_0a=x_0b$ for some $x_0\in S$ implies $xa=xb$ for all $x\in S$.

\begin{lemma}\label{lm3}(\cite{Nagy:sg-8}) A semigroup $S$ is left equalizer simple if and only if the right regular representation $S/\theta$ of $S$ is left cancellative.
\end{lemma}

\medskip

The previous lemmas show connections between $S$ and $S/\theta$, in special cases. In this paper we would like to find a common connection between the semigroups $S$, $S/\theta$ and $S/\theta ^*=(S/\theta)/(\theta ^*/\theta)$ in a general case. In our examination, the following construction plays an important role.

Let $S$ be a semigroup, $\Lambda$ an arbitrary set, and $P$ is a mapping of $\Lambda$ into $S$. It is easy to see that the set
$S\times \Lambda$ together with the operation $(s, \lambda )\circ _P (t, \mu )=(sP(\lambda )t, \mu )$ is a semigroup. This semigroup is the dual of the semigroup constructed in Exercise 6 for \S 8.2 of \cite{Clifford2:sg-20}. This semigroup will be denoted by $(S, \Lambda ; \circ _P)$.

In Section 2, we deal with the semigroups $(S, \Lambda ; \circ _P)$. We show that the left cancellativity and the right simplicity of a semigroup $S$ are inherited from $S$ to the semigroup $(S, \Lambda ; \circ _P)$.

In Section 3, the semigroups $(S, S/\theta ; \circ _P)$ are in the focus of our examination, where $S$ is an arbitrary semigroup and $P$ is an arbitrary mapping of the factor semigroup $S/\theta$ into $S$ with condition $P([s]_{\theta})\in [s]_{\theta}$. We show that, for an arbitrary semigroup $S$, the right regular representation of the semigroup $(S, S/\theta ; \circ _P)$ is isomorphic to the semigroup $(S/\theta ^*, S/\theta ; \circ _{P'})$, where $P'$ is the mapping of $S/\theta$ into $S/\theta ^*$ defined by $P'([s]_{\theta})= [s]_{\theta ^*}$.

In Section 4, we prove an embedding theorem for the semigroups $(S, S/\theta ; \circ _P)$, where $S$ is an idempotent-free left equalizer simple semigroup. We prove that if $S$ is a left equalizer simple semigroup without idempotent then the semigroup $(S, S/\theta ; \circ _P)$ can be embedded into a simple semigroup $(S", S/\theta ; \circ _{P"})$ containing a minimal left ideal.

For notations and notions not defined here, we refer to \cite{Clifford1:sg-1}, \cite{Clifford2:sg-20} and \cite{Nagybook:sg-5}.

\section{Hereditary properties}

In this section we show that the left cancellativity and the right simplicity of a semigroup $S$ are inherited from $S$ to the semigroup
$(S, \Lambda ; \circ _P)$.

\begin{lemma}\label{lm4} If $S$ is a left cancellative semigroup then the semigroup $(S, \Lambda , \circ _P)$  is also left cancellative
for any set $\Lambda$ and any mapping $P : \Lambda \mapsto S$.
\end{lemma}

\noindent
{\bf Proof}. Assume \[(s, \lambda )\circ _P(t, \mu )=(s, \lambda )\circ _P(r, \tau )\] for some $s, t, r\in S$ and $\lambda, \mu , \tau \in S$. Then
\[(sP(\lambda )t, \mu )=(sP(\lambda )r, \tau )\] from which it follows that \[sP(\lambda )t=sP(\lambda )r \quad \hbox{and}\quad \mu =\tau \] As $S$ is left cancellative, we get $t=r$. Hence \[(t, \mu )=(r, \tau ).\] Thus the semigroup $(S, \Lambda , \circ _P)$ is left cancellative.\hfill\openbox

\begin{lemma}\label{lm5} If $S$ is a right simple semigroup then the semigroup $(S; \Lambda ; \circ _P)$ is also right simple
for any set $\Lambda$ and any mapping $P : \Lambda \mapsto S$.
\end{lemma}

\noindent
{\bf Proof}. Let $(s, \lambda), (t, \mu)\in (S, \Lambda , \circ _P)$ be arbitrary elements. As $S$ is right simple, \[sP(\lambda )S=S.\] Then there is an element $x \in S$ such that \[sP(\lambda)x = t,\] and so \[(s, \lambda )\circ _P(x, \mu)=(sP(\lambda)x, \mu)=(t, \mu).\] From this it follows that the semigroup $(S, \Lambda ,\circ _P)$ is right simple.
\hfill\openbox

\begin{corollary}\label{cr1} If $S$ is a right group then the semigroup $(S; \Lambda ; \circ _P)$ is also a right group
for any semigroup $S$ and any mapping $P : \Lambda \mapsto S$.
\end{corollary}

\noindent
{\bf Proof}. As a semigroup is a right group if and only if it is right simple and left cancellative, our assertion follows from Lemma~\ref{lm4} and Lemma~\ref{lm5}.\hfill\openbox

\section{The right regular representation}

\bigskip

In this section we deal with the right regular representation of semigroups $(S, \Lambda ,\circ _P)$ in that case when $S$ is an arbitrary semigroup, $\Lambda$ is the factor semigroup $S/\theta$ and $P$ is an arbitrary mapping of $S/\theta$ into $S$ with condition that $P([s]_{\theta})\in [s]_{\theta}$ for every $s\in S$. We note that the product $\circ _P$ in the semigroup $(S, S/\theta , \circ _P)$ does not depend on choosing $P$, because $(s, [a]_{\theta})\circ _P(t, [b]_{\theta})=(sP([a]_{\theta})t, [b]_{\theta})$ for every $s, t, a, b\in S$, and $sa'=sa"$ for every $a', a"\in [a]_{\theta}$.

\begin{theorem}\label{th1} Let $S$ be an arbitrary semigroup. Let $P$ be a mapping of $S/\theta$ into $S$ with condition $P([a]_{\theta})\in [a]_{\theta}$ for every $[a]_{\theta}\in S/\theta$. Let $P'$ denote the mapping of $S/\theta$ onto $S/\theta ^*$ defined by $P'([a]_{\theta })=[a]_{\theta ^*}$. Then the right regular representation of the semigroup $(S, S/\theta ; \circ _P)$ is isomorphic to the semigroup
$(S/\theta ^*, S/\theta ; \circ _{P'})$.
\end{theorem}

\noindent
{\bf Proof}. Let $\theta ^{\bullet}$ denote the kernel of the right regular representation of the semigroup $(S, S/\theta , \circ _P)$. Let $\phi$ be the mapping of the factor semigroup $(S, S/\theta, \circ _P)/\theta ^{\bullet}$ onto the semigroup
$(S/\theta ^*, S/\theta, \circ _{P'})$ defined by
\[\phi ([(a, [b]_{\theta})]_{\theta ^{\bullet}})=([a]_{\theta ^*}, [b]_{\theta}),\]
where $[(a, [b]_{\theta})]_{\theta ^{\bullet}}$ denotes the $\theta ^{\bullet}$-class of $(S, S/\theta , \circ _P)$ containing the element
$(a, [b]_{\theta})$ of $(S, S/\theta, \circ _P)$.
To show that $\phi$ is injective, assume
\[\phi ([(a, [b]_{\theta})]_{\theta ^{\bullet}})=\phi ([(c, [d]_{\theta})]_{\theta ^{\bullet}})\]
for some $[(a, [b]_{\theta})]_{\theta ^{\bullet}}, [(c, [d]_{\theta})]_{\theta ^{\bullet}}\in (S, S/\theta, \circ _P)/\theta ^{\bullet}$.
Then
\[([a]_{\theta ^*}, [b]_{\theta})=([c]_{\theta ^*}, [d]_{\theta})\] and so

\[[a]_{\theta ^*}=[c]_{\theta ^*}\quad \hbox{and}\quad [b]_{\theta}=[d]_{\theta}.\]

Let $x, \xi \in S$ be arbitrary elements. Then $x\xi a=x\xi c$ and so
\[(x, [\xi ]_{\theta})(a, [b]_{\theta})=(x\xi a, [b]_{\theta})=(x\xi c, [d]_{\theta})=(x, [\xi ]_{\theta})(c, [d]_{\theta}).\]
Hence
\[((a, [b]_{\theta}), (c, [d]_{\theta}))\in \theta ^{\bullet},\] that is,
\[ [(a, [b]_{\theta})]_{\theta ^{\bullet}}=[(c, [d]_{\theta})]_{\theta ^{\bullet}}.\]
Thus $\phi$ is injective. Consequently $\phi$ is bijective.
It remains to show that $\phi$ is a homomorphism. Let
\[[(a, [b]_{\theta})]_{\theta ^{\bullet}}, [(c, [d]_{\theta})]_{\theta ^{\bullet}}\in (S, S/\theta, \circ _P)/\theta ^{\bullet}\] be arbitrary. Then \[\phi ([(a, [b]_{\theta})]_{\theta ^{\bullet}}[(c, [d]_{\theta})]_{\theta ^{\bullet}})=\phi ([(abc, [d]_{\theta})]_{\theta ^{\bullet}})=\]
\[=([abc]_{\theta ^*}, [d]_{\theta})=([a]_{\theta ^*}[b]_{\theta ^*}[c]_{\theta ^*}, [d]_{\theta})=\]
\[=([a]_{\theta ^*}, [b]_{\theta})\circ _{P'}([c]_{\theta ^*}, [d]_{\theta})=\phi ([(a, [b]_{\theta})]_{\theta ^{\bullet}})\circ _{P'}\phi ([(a, [b]_{\theta})]_{\theta ^{\bullet}}).\]
Hence $\phi$ is a homomorphism, and so it is an isomorphism of the right regular representation of the semigroup $(S, S/\theta, \circ _P)$ onto the semigroup $(S/\theta ^*, S/\theta, \circ _{P'})$.\hfill\openbox

\begin{corollary}\label{cr2} If $S$ is a left group then the semigroup $(S, S/\theta, \circ _P)$ is M-inversive.
\end{corollary}

\noindent
{\bf Proof}. If $S$ is a left group then $S/\theta$ is a group by Lemma~\ref{lm1}. As $S^2=S$, we have $\theta=\theta ^*$. Thus $S/\theta ^*$ is a group. By Corollary~\ref{cr1}, the semigroup
$(S/\theta ^* , S/\theta ; \circ _{P'})$ is a right group and so, by Theorem~\ref{th1} and Lemma~\ref{lm2}, $(S, S/\theta , \circ _P)$ is an M-inversive semigroup.
\hfill\openbox

\begin{corollary}\label{cr3} If $S$ is an M-inversive semigroup then the semigroup $(S, S/\theta, \circ _P)$ is M-inversive.
\end{corollary}

\noindent
{\bf Proof}. If $S$ is M-inversive then $S/\theta$ is a right group by Lemma~\ref{lm2}. As the kernel of the right regular representation of a right group is the identity relation, $\theta ^*/\theta$ is the identity relation on $S/\theta$. Then $\theta ^*=\theta$ and so $S/\theta ^*$ is a right group. By Corollary~\ref{cr1}, $(S/\theta ^* \times S/\theta ; \circ _{P'})$ is a right group. By Theorem~\ref{th1} and Lemma~\ref{lm2},
$(S, S/\theta , \circ _P)$ is an M-inversive semigroup.\hfill\openbox

\begin{corollary}\label{cr4} If $S$ is a left equalizer simple semigroup then the semigroup $(S, S/\theta, \circ _P)$
is left equalizer simple.
\end{corollary}

\noindent
{\bf Proof}. Let $S$ be a left equalizer simple semigroup. Then, by Lemma~\ref{lm3}, $S/\theta$ is a left cancellative semigroup. As the right regular representation of a left cancellative semigroup is the identity relation, $\theta ^*/\theta$ is the identity relation on $S/\theta$. From this it follows that $\theta ^*=\theta$ and so $S/\theta ^*$ is a left cancellative semigroup. By Lemma~\ref{lm4}, $(S/\theta ^*, S/\theta, \circ _{P'})$ is a left cancellative semigroup. By Theorem~\ref{th1} and by Lemma~\ref{lm3}, the semigroup $(S, S/\theta , \circ _P)$ is M-inversive.\hfill\openbox

\section{An embedding theorem}

\bigskip

In this section we deal with such semigroups $(S, S/\theta , \circ _P)$ in which $S$ is an idempotent-free left equalizer simple semigroup.

\begin{theorem}\label{th2} If $S$ is a left equalizer simple semigroup without idempotents then the semigroup $(S, S/\theta, \circ _P)$ can be embedded into a simple semigroup $(S", S/\theta, \circ _{P"})$ containing a minimal left ideal.
\end{theorem}

\noindent
{\bf Proof}. Let $S$ be a left equalizer simple semigroup without idempotent. Then, by Theorem 8.19 of \cite{Clifford2:sg-20}, there is an embedding $\tau$ of $S$ into a left simple semigroup $S"$ without idempotents. Consider the semigroup $(S", S/\theta , \circ _{P"})$, where the mapping $P" : S/\theta \mapsto S"$ is defined by
$P"([s]_{\theta})=\tau (P[s]_{\theta})$. By the dual of Exercise 6 for \S 8.2 of \cite{Clifford2:sg-20},
$(S", S/\theta , \circ _{P"})$ is a simple semigroup containing a minimal left ideal. We show that
\[\Phi : (a, [b]_{\theta}) \mapsto (\tau (a), [b]_{\theta})\] is an embedding of the semigroup
$(S, S/\theta, \circ _P)$ into the semigroup $(S", S/\theta , \circ _{P"})$.

First we show that $\Phi$ is injective. Assume
\[ \Phi ((a, [b]_{\theta}))=\Phi ((c, [d]_{\theta}))\] for some $a, b, c, d\in S$. Then
\[\tau (a)=\tau (c)\quad \hbox{and}\quad [b]_{\theta}=[d]_{\theta}.\]
As $\tau$ is injective, we get \[(a, [b]_{\theta})=(c, [d]_{\theta}).\]

Next we show that $\Phi$ is a homomorphism. Let \[(a, [b]_{\theta}), (c, [d]_{\theta}) \in (S, S/\theta , \circ _P)\] be arbitrary elements.
Then
\[\Phi ((a, [b]_{\theta})\circ _P (c, [d]_{\theta}))=\Phi ((aP([b]_{\theta})c, [d]_{\theta}))=(\tau (aP([b]_{\theta})c), [d]_{\theta})=\]
\[(\tau (a)\tau (P([b]_{\theta}))\tau (c), [d]_{\theta })=(\tau (a)P"([b]_{\theta})\tau (c), [d]_{\theta })=\]
\[=(\tau (a) ,[b]_{\theta})\circ _{P"} (\tau (c), [d]_{\theta})=\Phi((a, [b]_{\theta}))\circ _{P"} \Phi ((c, [d]_{\theta}))\]
and so $\Phi$ is a homomorphism. Consequently $\Phi$ is an embedding of the semigroup
$(S, S/\theta, \circ _P)$ into the simple semigroup $(S", S/\theta , \circ _{P"})$ containing a minimal left ideal.\hfill\openbox

\end{document}